\documentclass[11pt]{article}

\newfont{\bb}{msbm10 at 11pt}
\def\r{\hbox{\bb R}}
\def\h{\hbox{\bb H}}
\def\s{\hbox{\bb S}}

\newenvironment{proof}{\trivlist
\item[\hskip\labelsep{\it Proof}\,:]}{\hfill{$q.e.d.$}\endtrivlist}
\newtheorem{theorem}{Theorem}[section]
\newtheorem{lemma}[theorem]{Lemma}

\newtheorem{corollary}[theorem]{Corollary}

\newtheorem{remark}[theorem]{Remark}

\usepackage{graphicx}

\begin{document}

\title{Parabolic Weingarten surfaces in hyperbolic space}
\author{Rafael L\'opez\footnote{Partially
supported by MEC-FEDER
 grant no. MTM2007-61775.}\\
Departamento de Geometr\'{\i}a y Topolog\'{\i}a\\
Universidad de Granada\\
18071 Granada (Spain)\\
e-mail:\ {\tt rcamino@ugr.es}}
\date{}
\maketitle

\begin{abstract}
A surface in hyperbolic space $\h^3$ invariant by a group of
parabolic isometries is called a parabolic surface. In this paper we
investigate parabolic surfaces of $\h^3$ that satisfy a linear
Weingarten  relation of the form $a\kappa_1+b\kappa_2=c$ or
$aH+bK=c$, where $a,b,c\in \r$ and, as usual, $\kappa_i$ are the
principal curvatures, $H$ is the mean curvature and $K$ is de
Gaussian curvature. We classify all parabolic  linear  Weingarten
surfaces in hyperbolic space.
\end{abstract}


\section{Introduction}\label{intro}

A surface $S$ in 3-dimensional hyperbolic space $\h^3$ is called a
{\it Weingarten surface} if there is some relation between its two
principal curvatures $\kappa_1$ and $\kappa_2$, that is, if there is
a smooth function $W$ of two variables such that
$W(\kappa_1,\kappa_2)=0$. In particular, if $K$ and $H$  denote
respectively  the Gauss curvature and the mean curvature of  $S$,
the identity $W(\kappa_1,\kappa_2)=0$  implies a relation
$U(K,H)=0$.  In this paper we study Weingarten surfaces that satisfy
the simplest case for $W$ and $U$, that is,  of linear type:
\begin{equation}\label{wein1}
a\ \kappa_1+b\ \kappa_2=c
\end{equation}
and
\begin{equation}\label{wein2}
a\ H+b\ K=c,
\end{equation}
where $a,b,c\in\r$.  We say in both cases that $S$ is a {\it linear
Weingarten surface} and we abbreviate  by $LW$-surface. In the set
of $LW$-surfaces, it is worth mentioning three families of surfaces
that correspond with trivial choices of  $a,b$ and $c$:
\begin{enumerate}
\item Umbilical surfaces, when  $a=-b$ and $c=0$ in (\ref{wein1}).
\item Surfaces with constant mean curvature: they appear if we choose $a=b$ in (\ref{wein1}) or $b=0$ in (\ref{wein2}).
\item Surfaces with constant Gaussian curvature, with
the choice  $a=0$ in (\ref{wein2}).
\end{enumerate}
We call these three families of surfaces as {\it trivial
$LW$-surfaces}. Although these three kinds of surfaces have been
studied in the literature, the  classification of $LW$-surfaces in
the general case  is almost completely open today. One of the
objectives of this work is to provide examples of new surfaces. The
idea is to consider  surfaces of revolution since in such case, the
relations (\ref{wein1}) and (\ref{wein2}) reduce to an ordinary
differential equation that describes the shape of the profile curve
that generates the surface. In the literature, part of this work  has been carried in \cite{cor,gmm,st3,yo}.

In hyperbolic space $\h^3$,  there exist three kinds of rotational
surfaces depending on the type of orbits under the action of the
group of isometries: spherical, hyperbolic and parabolic rotational
surfaces. Our interest are the  surfaces invariant by a group of
parabolic isometries. A parabolic group of isometries of $\h^3$ is
formed by isometries that leave fix one double point of the ideal
boundary $\s^2_{\infty}$  of $\h^3$. We say that a  surface is a
{\it parabolic surface} of $\h^3$ if it is invariant by a group of
parabolic isometries.  A such parabolic surface $S$ is determined by
a generating curve  $\alpha$ obtained  by the intersection of $S$
with any geodesic plane orthogonal to the orbits of the group.
Parabolic surfaces in $\h^3$  were introduced by Do Carmo and
Dajczer in \cite{dd} focusing in the study of surfaces with constant
mean curvature (see also, \cite{go}). More exactly, and with respect
to the trivial $LW$-surfaces, we point out that umbilical surfaces
in $\h^3$ are well known (see for example \cite{sp});  parabolic
surfaces with constant mean curvature are given in the cited  papers
\cite{dd} and \cite{go}, and finally, parabolic surfaces with
constant Gaussian curvature are described
 in   \cite{lo2,mo}.

Between the facts  interesting for LW-surfaces of
parabolic type, we point out some of them. First, the  question
whether the surface can be extended to be complete, which it is
given in terms of the generating curve. Second, if a complete
parabolic Weingarten surface is embedded. For example, this occurs
if the surface has constant Gaussian curvature \cite{lo2,mo}. However, there exist constant mean curvature
 non-embedded surfaces that are complete  \cite{go}.
Finally, the question about the behavior of the surface in relation
with the ideal boundary $\s^2_{\infty}$. We know that the
asymptotic boundary of surface contains the fixed point of the
parabolic group of isometries.

This paper is organized as follows. In Section \ref{spre} we
establish the differential equations that govern the  parabolic LW-surfaces and some properties about their
symmetries. In Sections \ref{sw1} and \ref{sw2} we study all
parabolic surfaces in $\h^3$ that satisfy equations (\ref{wein1})
and (\ref{wein2}), respectively. We give a complete description of
such surfaces, which depends on certain relations of the parameters
$a,b$ and $c$. In all the cases, we assume that the
generating curve of the parabolic surface has a tangent line
parallel to the ideal boundary $\s^2_{\infty}$. In Section
\ref{final}, we point out how this assumption can be dropped, which
will complete the classification.

For the explicit classification, we refer the readers to Sections
\ref{sw1}, \ref{sw2} and \ref{final}. However, we can announce  some facts that are worth to point out.

\begin{quote} {\it Any parabolic surface in $\h^3$ that satisfies the
relation $a\kappa_1+b\kappa_2=c$ can be extended to be  complete.
The asymptotic boundary of any such surface is one point, one circle
 of two tangent circles. If it  is one circle,  the surface
 is umbilical. Moreover, there exist  surfaces that are graphs on
   $\s^2_{\infty}$.}
 \end{quote}

  \begin{quote}{\it There exist complete parabolic surfaces in $\h^3$
 that satisfy $aH+bK=c$. For these surfaces, the asymptotic boundary
 is the point $\infty$, one circle or two circles tangent at $\infty$. Some of the above
 surfaces are graphs on  $\s^2_{\infty}$. There exist surfaces that can not extend
 to be complete}
  \end{quote}

\section{Preliminaries and first properties}\label{spre}

In this section we fix some notations and we give some properties about
the symmetries of parabolic $LW$-surfaces.
Let us consider the upper half-space model of the
hyperbolic three-space $\h^3$, namely,
$${\h}^{3}=:{\r}^3_+=\{(x,y,z)\in{\r}^{3};z>0\}$$
equipped with the metric
$$\langle,\rangle=\frac{dx^{2}+dy^{2}+dz^{2}}{z^{2}}.$$
In what follows, we will  use the words  "vertical" or "horizontal" in the
usual affine sense of $\r^3_+$.
The ideal boundary $\s^2_{\infty}$ of $\h^3$ is identified with the one point compactification of the plane $\Pi\equiv\{z=0\}$, that is,
$\s^2_{\infty}=\Pi\cup\{\infty\}$ and it  corresponds with the asymptotic classes of geodesics rays of $\h^3$. The asymptotic boundary of a set
$\Sigma\subset\h^3$ is defined as
$\partial_{\infty}\Sigma=\overline{\Sigma}\cap \s^2_{\infty}$,
where $\overline{\Sigma}$ is the closure of $\Sigma$ in $\{z\geq 0\}\cup\{\infty\}$. Let $L=\{(x,0,0),x\in\r\}$.

A parabolic group  of isometries $G$ of $\h^3$ is a group of isometries that
admits a fixed double point at $\s^2_{\infty}$.
These isometries leave globally fixed each horocycle tangent to the fixed point.
 In our model, and without loss of generality, we  take the point  $\infty$ of
$\s^2_{\infty}$ as the point that fixes $G$.
 Then the group $G$ is defined by the horizontal (Euclidean) translations in the
direction of a  horizontal vector $\xi$ with $\xi\in\Pi$:
$G=\{T_a; a\in\r, T_a(p)=p+a\xi\}$. The orbits are then horizontal
straight lines parallel to $\xi$. We can also view this group as the set of
reflections with respect to any geodesic plane orthogonal to $\xi$. Actually, the
parabolic group $G$ is generated by all reflections with respect to the geodesic
planes  orthogonal to $\xi$. The space of orbits is then represented in any geodesic plane
of this family. This will be done in our study.

Let $G$ be a group of parabolic isometries. Without loss of
generality, we assume that the horizontal vector $\xi$ that defines
the group of is the vector $\xi=(0,1,0)$. Let $P=\{(x,0,z);z>0\}$,
which it is a vertical geodesic plane  orthogonal to $\xi$. Then a
surface $S$ invariant by $G$ intersects $P$ in a curve $\alpha$
called the {\it generating curve} of $S$. If $S$ is a parabolic
$LW$-surface, we shall obtain an ordinary differential equation for
the curve $\alpha$, equations (\ref{wein11}) and (\ref{wein22})
below. If we assume that $S$ is a complete surface, the
possibilities about its asymptotic boundary $\partial_\infty S$ are:
a circle ($\partial_{\infty}\alpha$ is a point of $L$ or one point
of $L$  together $\infty$), two tangent circles
($\partial_{\infty}\alpha$ are two different points of $L$)
 or it is one point ($\partial_{\infty}\alpha=\emptyset$ or $\infty$).

Let $S$ be a parabolic (connected) surface in $\h^3$ and let
$X(s,t)=(x(s),t,z(s))$ be a parametrization of $S$, where $t\in\r$
and the curve $\alpha$ will be assumed to be parametrized by the arc
length with respect to the Euclidean metric, whose domain of
definition $I$ is an open interval of real numbers including zero.
The principal directions at each point are $\partial_s X$ and
$\partial_t X$. Denote $\theta$ the angle that makes the velocity
$\alpha'(s)$ with the $x$-axis, that is, $x'(s)=\cos\theta(s)$ and
$z'(s)=\sin\theta(s)$ for a certain differentiable function
$\theta$. The derivative $\theta'(s)$ of the function $\theta(s)$ is
the Euclidean curvature of $\alpha$. From the hyperbolic viewpoint,
the hyperbolic curvature of $\alpha$ at $s$ is exactly
$z(s)\theta'(s)+\cos\theta(s)$.

Consider the Gauss map
$N(s,t)$ induced by the immersion $X(s,t)$, that is, $N(s,t)=z(s)(-\sin\theta(s),0,\cos\theta(s))$. Then  the principal curvatures $\kappa_i$
of $S$ are
\begin{equation}\label{k1k2}
\kappa_1(s,t)=z(s)\theta'(s)+\cos\theta(s),\hspace*{1cm}\kappa_2(s,t)=\cos\theta(s)
\end{equation}
and  the mean curvature $H=\frac{\kappa_1+\kappa_2}{2}$ and Gaussian curvature $K=\kappa_1\kappa_2-1$ are
\begin{equation}\label{hk}
H(s,t)=\frac{z(s)}{2}\theta'(s)+\cos\theta(s),\hspace*{1cm}
K(s,t)=z(s)\cos\theta(s)\theta'(s)-\sin\theta(s)^2.
\end{equation}
Thus, parabolic $LW$-surfaces in $\h^3$ are given by  curves $\alpha$ whose coordinate functions satisfy
\begin{equation}\label{alpha}
 \left\{\begin{array}{lll}
 x'(s)&=& \displaystyle \cos\theta(s)\\
 z'(s)&=&\displaystyle \sin\theta(s)
\end{array}
\right.
\end{equation}
together the equation
\begin{equation}\label{wein11}
a z(s)\theta'(s)+(a+b)\cos\theta(s)=c
\end{equation} or
\begin{equation}\label{wein22}
\left(\frac{a}{2}+b\cos\theta(s)\right)z(s)\theta'(s)+a\cos\theta(s)-b\sin^2\theta(s)=c
\end{equation}
depending  if $S$ satisfies the Weingarten relation
 (\ref{wein1}) or (\ref{wein2}) respectively.
We consider the initial conditions
\begin{equation}\label{eq2}
x(0)=0,\hspace*{.5cm}z(0)=z_0>0,\hspace*{.5cm}\theta(0)=\theta_0.
\end{equation}

We first prove two properties about the symmetries of the solutions of (\ref{wein11}) and
 (\ref{wein22}).

\begin{lemma}\label{si1}
Let $\alpha$ be a solution of the initial value problem (\ref{alpha})-(\ref{wein11})
or (\ref{alpha})-(\ref{wein22}). Suppose that
$z'(s_0)=0$ for a real number $s_0$. Then $\alpha$ is symmetric with respect to
the vertical line $x=x(s_0)$ of the $xz$-plane.
\end{lemma}

\begin{proof}
We do the proof for a solution of (\ref{alpha})-(\ref{wein11}) and
the reasoning is analogous in the another case.
Since $\sin\theta(s_0)=0$, then $\theta(s_0)=k\pi$ for some integer number $k$.
The triplets of functions $\{x(s_0+s),z(s_0+s),\theta(s_0+s)\}$ and
$\{2x(s_0)-x(s_0-s),z(s_0-s),-\theta(s_0-s)+2k\pi\}$ satisfy the same differential equations and  the same initial conditions at $s=0$. The uniqueness of solutions concludes the result.

\end{proof}

\begin{lemma}\label{cero}
Let $\alpha$ be  a solution of the initial value problem (\ref{alpha})-(\ref{wein11})
or (\ref{alpha})-(\ref{wein22}). Suppose that $\theta'(s_0)=0$
for a real number $s_0$. Then $\alpha$ is a straight line and the corresponding surface is a totally geodesic plane, an
equidistant surface or a horosphere.
\end{lemma}

\begin{proof} As in  Lemma \ref{si1}, we restrict to the case that $\alpha$ satisfies
(\ref{alpha})-(\ref{wein11}). If $\{x(s),z(s),\theta(s)\}$ is a such
solution, then
$$\{\cos\theta(s_0)(s-s_0)+x(s_0),\sin\theta(s_0)(s-s_0)+z(s_0),\theta(s_0)\}$$
is a solution of (\ref{alpha})-(\ref{wein11}) with the same initial
conditions at $s=s_0$. Thus these three functions are the very
solutions of the differential equations system.
\end{proof}

Finally, and to end with this section, we consider the relation $a\kappa_1+b\kappa_2=c$ in the
 case that $a$ or $b$ is zero. Then the surface has
one constant principal curvature. Actually, each orbit $t\mapsto X(s,t)$ is
a line of curvature, whose curvature, namely $\cos\theta(s)$, is constant  along the line.
On the other hand, the normal curvature of the curve  $s\mapsto X(s,t)$ agrees with the
(hyperbolic) curvature as planar curve in $\h^3$. Thus, if it is constant, it is
well known that the curve is a straight line or a Euclidean circle. We can see this as follows.

\begin{theorem}\label{principal}
 The only parabolic surfaces in $\h^3$ with one constant principal curvature
are totally geodesic planes, equidistant surfaces, horospheres and
Euclidean horizontal right-cylinders.
\end{theorem}

\begin{proof} We distinguish between the two principal curvatures $\kappa_1$ and $\kappa_2$.
Assume that $\kappa_1=c$, where $c$ is a constant. Then $\theta'(s)
z(s)=c-\cos\theta(s)$. By differentiation of this expression and
using (\ref{alpha}) we obtain $\theta''(s)=0$ for all $s$. Then
$\theta'$ is constant and hence that  from the Euclidean viewpoint,
the curve is a piece of a straight line or a circle, which generates
(pieces of) geodesic planes, equidistant surfaces, horospheres and
horizontal right-cylinders.

Suppose now that $\kappa_2$ is constant, that is, $\cos\theta(s)=c$. This means that
$\theta$ is constant, and so, $\alpha$ is a straight line. This gives
totally geodesic planes (if $c=0$), equidistant surfaces (if $0<|c|<1$) and
horospheres (if $|c|=1$).
\end{proof}

After an isometry of the ambient space, the surfaces that are
Euclidean horizontal right-cylinders are banana-shaped surfaces
whose end points agree at one point of $\s^2_{\infty}$.

\section{Parabolic surfaces satisfying $ \kappa_1=m\kappa_2+n$.}\label{sw1}

In this section we shall consider parabolic surfaces that satisfy
the relation (\ref{wein1}). The case that  one of the principal
curvatures $\kappa_i$ is constant   has been completely studied in
Theorem \ref{principal}. Thus we deal with the case that both $a$ and
$b$ are non-zero numbers. Then the relation (\ref{wein1}) can written as
\begin{equation}\label{w1}
\kappa_1=m\kappa_2+n
\end{equation}
where $m,n\in\r$, $m\not=0$. By using (\ref{k1k2}), we have
\begin{equation}\label{w11}
\theta'(s)=\frac{(m-1)\cos\theta(s)+n}{z(s)}.
\end{equation}
After a change of orientation on the surface, we suppose that $n\geq 0$. We discard the trivial $LW$-surfaces, that is,
umbilical surfaces corresponding to $(m,n)=(1,0)$ and the surfaces
with constant mean curvature, that is,  $m=-1$. We consider
$\theta(0)=\theta_0=0$ in the  initial condition. In particular and
from Lemma \ref{si1}, the generating curve $\alpha$ is symmetric
with respect to the line $x=0$ of the $xz$-plane $P$. Multiplying in
(\ref{w11}) by $\sin\theta$ and integrating, we obtain
\begin{equation}\label{integral2}
n+\cos\theta(s)=\frac{2-m}{z(s)}\int_0^s
\left(\sin\theta(t)\cos\theta(t)\right)\ dt+(n+1) \frac{z_0}{z(s)}.
\end{equation}
Equation (\ref{w11}) yields at $s=0$,
$$\theta'(0)=\frac{n+m-1}{z_0}.$$
By Lemma \ref{cero}, if the function $\theta'(s)$ vanishes at some point $s$,
then $\theta'=0$ and $\alpha$ is a straight line. If
$\theta'(0)\not=0$, then $\theta(s)$ is a strictly monotonic
function on $s$. Let $(-\bar{s},\bar{s})$ be the maximal domain of
solutions of (\ref{alpha})-(\ref{w11}) under the initial conditions
(\ref{eq2}). Denote $\theta_1=\lim_{s\rightarrow\bar{s}}\theta(s)$.
Depending on the sign of $\theta'(0)$, we consider three cases.

\subsection{Case $n+m-1>0$}

Here $\theta'(0)>0$ and so, $\theta$ is strictly increasing in its
domain.

\begin{enumerate}
\item Subcase $m<n+1$. In particular, $n>0$.
 We prove that $\theta$ attains the value $\pi/2$.
Assume on the contrary, that is, $\theta_1\leq\pi/2$ and we will
arrive to a contradiction. As $z'(s)=\sin\theta(s)>0$,  $z(s)$ is
strictly increasing in $(0,\bar{s})$. Then $z(s)\geq z_0$ and  the
derivatives of $\{x(s),z(s),\theta(s)\}$ in equations
(\ref{alpha})-(\ref{w11}) are bounded. This means that
$\bar{s}=\infty$. As
$\lim_{s\rightarrow\infty}z'(s)=\sin\theta_1>0$, then
$\lim_{s\rightarrow\infty}z(s)=\infty$. Let $s\rightarrow\infty$ in
(\ref{integral2}). If the integral that appears in the right-side is
bounded, then $n+\cos\theta_1=0$, that is, $\cos\theta_1=n=0$:
contradiction. If the integral is not bounded, and using the
L'H\^{o}pital's rule, $n+\cos\theta_1=(2-m)\cos\theta_1$, that is,
$(m-1)\cos\theta_1+n=0$. Then $m-1\leq 0$ and the hypothesis
$n+m-1>0$ yields $\cos\theta_1= n/(1-m)>1$: contradiction.

Therefore, there exists a first value $s_1$ such that
$\theta(s_1)=\pi/2$.  We prove that $\theta(s)$ attains the value
$\pi$. By contradiction, we assume
  $\theta_1\leq\pi$ and  $z(s)$
is strictly increasing again. We then have $\bar{s}=\infty$ again
and $\theta'(s)\rightarrow 0$ as $s\rightarrow\infty$. If $z(s)$ is
bounded, then (\ref{integral2}) implies $(m-1)\cos\theta_1+n=0$. As
$m-1=n=0$ is  impossible, then $m-1>0$ since $\cos\theta_1<0$. But
the hypothesis $m<n+1$ implies that $\cos\theta_1=-n/(m-1)<-1$,
which it is a contradiction. Thus $z(s)\rightarrow\infty$ as
$s\rightarrow\infty$. By using (\ref{integral2}) again, and letting
$s\rightarrow\infty$, we have $n+\cos\theta_1=0$. In particular,
$0<m<2$. For the contradiction, we obtain a second integral from
(\ref{w11}) multiplying by $\cos\theta(s)$:
$$\sin\theta(s)=\frac{s}{z(s)}+\frac{1}{z(s)}
\int_0^s\left( n\cos\theta(t)+(m-2)\cos^2\theta(t)\right)\ dt.$$ If
the integral is bounded, then $\sin^2\theta_1=1$: contradiction.
Thus, the integral  is not bounded and L'H\^{o}pital rule implies
$\sin^2\theta_1=1+n\cos\theta_1+(m-2)\cos^2\theta_1$. This equation,
together $n+\cos\theta_1=0$ yields $(m-2)\cos^2\theta_1=0$:
contradiction.

As conclusion,  there exists a first value $s_2$ such that
$\theta(s_2)=\pi$. By Lemma \ref{si1}, the curve $\alpha$ is
symmetric with respect to the line $x=x(s_2)$.  Moreover, and
putting $T=2s_2$, we have:
$$x(s+T)=x(s)+x(T),\hspace*{.5cm}z(s)=z(s+T),\hspace*{.5cm}
\theta(s+T)=\theta(s)+2\pi.$$ This means that $\alpha$ is invariant
by a group of horizontal translations orthogonal to the orbits of
the parabolic group.

\item Subcase $m\geq n+1$. With this hypothesis and as
$\theta'(s)>0$, the equation (\ref{w11})
 implies that $\cos\theta(s)\not=-1$ for any $s$. Thus $-\pi <\theta(s)<\pi$.
For $s>0$, $z'(s)=\sin\theta(s)>0$ and then $z(s)$ is increasing on
$s$ and so, $\theta'(s)$ is a bounded function. This implies
$\bar{s}=\infty$.
 We show that either there exists $s_0>0$ such $\theta(s_0)=\pi/2$ or
$\lim_{s\rightarrow\infty}\theta(s)=\pi/2$.

As in the above subcase, and with the same notation,  if
$\theta(s)<\pi/2$ for any $s$, then $n+\cos\theta_1=0$ or
$(m-1)\cos\theta_1+n=0$.  As $\cos\theta_1\geq 0$ and since $m-1\geq
n$, it implies that this occurs if and only if $n=0$ and
$\theta_1=\pi/2$. In such case, $z''(s)=\theta'(s)\cos\theta(s)>0$,
that is, $z(s)$ is a convex function.
 As conclusion, if
$n>0$, there exists a value  $s_0$ such that $\theta(s_0)=\pi/2$, and
there exists $\theta_1\in (\pi/2,\pi]$ such that
$\lim_{s\rightarrow\infty}\theta(s)=\theta_1$.
\end{enumerate}

\begin{figure}[htbp]\begin{center}
\includegraphics[width=6cm]{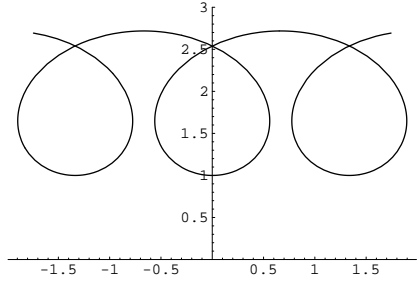}
\end{center}
\caption{The generating curves of a parabolic
surfaces with $\kappa_1=m\kappa_2+n$ with $n+m-1>0$ and subcase $m<n+1$. Here
$m=1$ and $n=2$. The value $z_0$ is $z_0=1$.}\label{11}
\end{figure}

\begin{figure}[htbp]\begin{center}
\includegraphics[width=6cm]{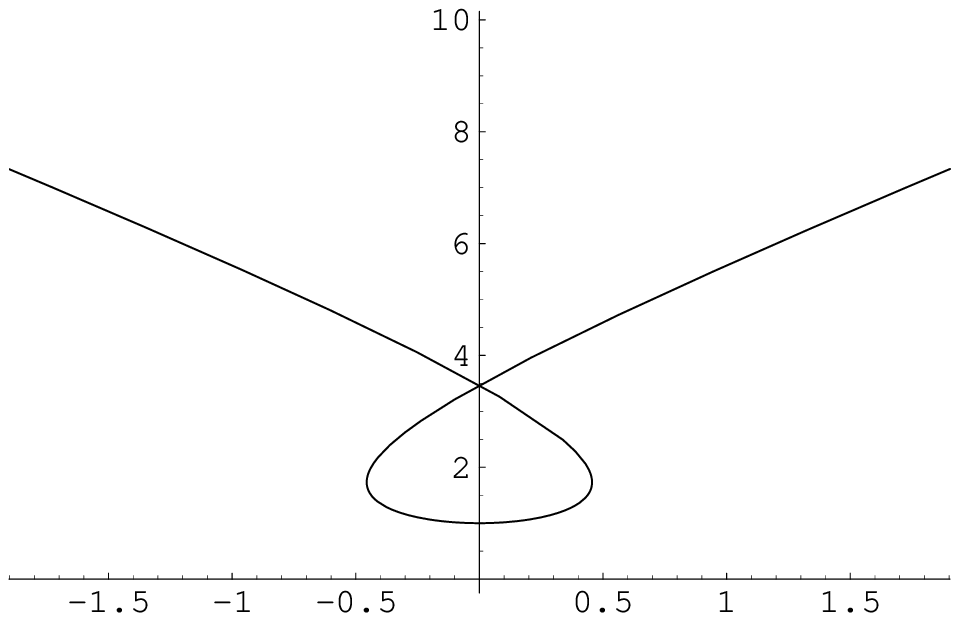}\hspace*{2cm}
\includegraphics[width=6cm]{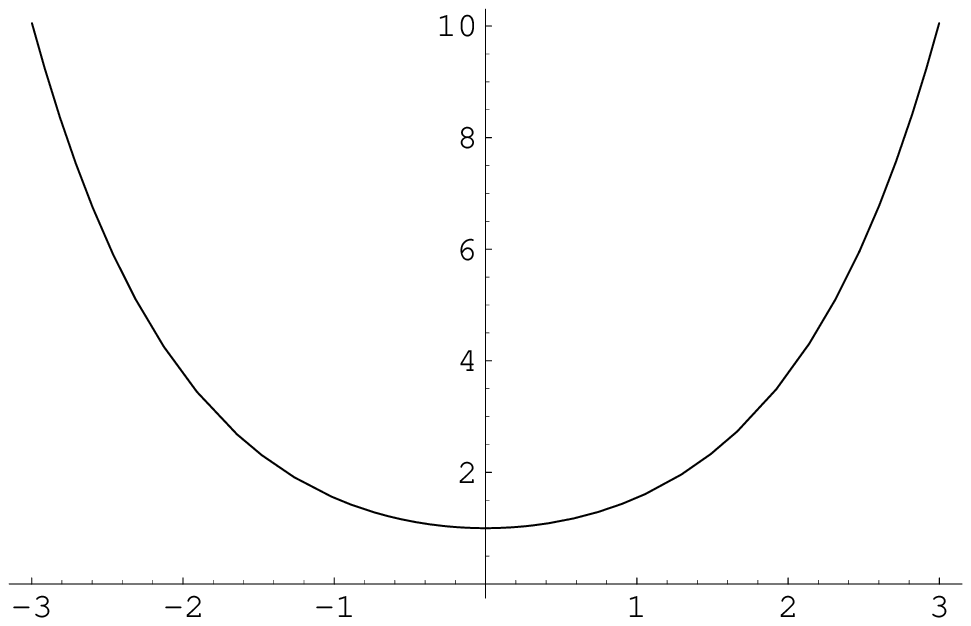}
\end{center}
\begin{center}(a)\hspace*{8cm}(b)\end{center}
\caption{The generating curves of a parabolic surfaces with
$\kappa_1=m\kappa_2+n$ with $n+m-1>0$ and subcase $m\geq n+1$.  In
(a), we have $m=3$ and $n=1$; in (b), $m=2$ and $n=0$. The value
$z_0$ is $z_0=1$.}\label{21}
\end{figure}

\begin{theorem}\label{th1} Let $\alpha(s)=(x(s),0,z(s))$ be the generating curve of a
LW-parabolic surface  $S$ in hyperbolic space $\h^3$ whose principal
curvatures  satisfy the relation $\kappa_1=m\kappa_2+n$. Consider $n\geq 0$ and
that  $\theta(0)=0$ in the initial condition (\ref{eq2}).
Assume $n+m-1>0$.
\begin{enumerate}
\item If $m<n+1$, then $\alpha$ is invariant by a group of translations in the
$x$-direction. Moreover,  $\alpha$ has self-intersections and it
presents one maximum and one minimum in each period, with vertical
points  between maximum and minimum. The velocity $\alpha'$ twirls
around the origin. See Fig. \ref{11}.
\item Assume $m\geq n+1$.
If $n>0$, then $\alpha$ has a minimum with self-intersections (see
Fig. \ref{21}, case (a)). If $n=0$, then $\alpha$ is a convex graph
on $L$, with a minimum (see Fig. \ref{21}, case (b)).
\end{enumerate}
\end{theorem}


\subsection{Case $n+m-1=0$}


In the case $n+m-1=0$, $\theta'(0)=0$ and thus, $\theta'(s)=0$ for
any $s$. As $\theta(0)=0$, then $\theta(s)=0$ for  any $s$. This
implies that $\alpha(s)$ is a horizontal straight line.

\begin{figure}[htbp]\begin{center}
\includegraphics[width=6cm]{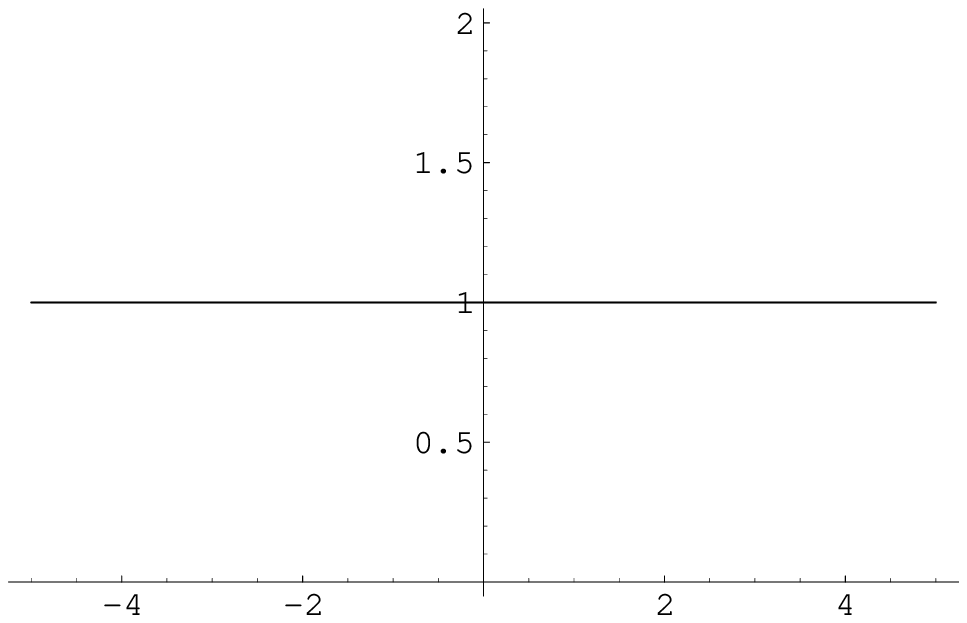}\hspace*{2cm}\includegraphics[width=6cm]{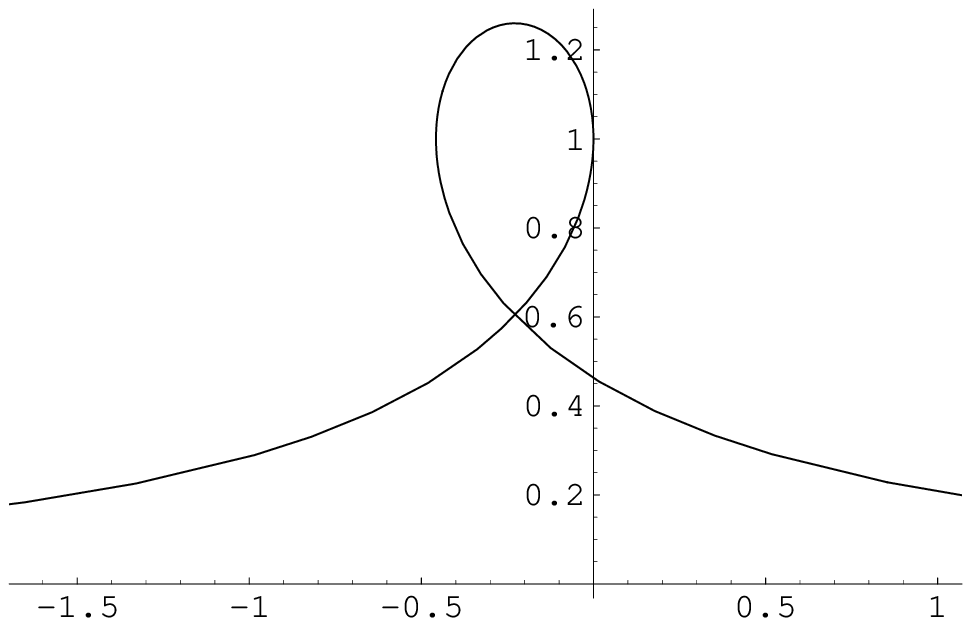}\end{center}
\begin{center}(a)\hspace*{8cm}(b)\end{center}
\caption{The generating curve of a parabolic surface with
$\kappa_1=m\kappa_2+n$ with $n+m-1=0$. Here $z_0=1$ and $m=-2$ and
$n=3$.  In the case (a), $\theta(0)=0$ and in the case (b),
$\theta(0)=\pi/2$ (see Theorem \ref{ultimo}).}\label{41}
\end{figure}

\begin{theorem} Let $\alpha(s)=(x(s),0,z(s))$ be the generating curve of a
parabolic surface  $S$ in hyperbolic space $\h^3$. Assume that the
principal curvatures of $S$ satisfy the relation
$\kappa_1=m\kappa_2+n$ with $n+m-1=0$ and $n\geq 0$. If
$\theta(0)=0$ in the initial condition (\ref{eq2}), then $S$ is  a
horosphere (see Fig. \ref{41}, case (a)).
\end{theorem}

\subsection{Case $n+m-1<0$}


With this assumption, $\theta(s)$ is a decreasing function. As
$n\geq 0$ and from (\ref{w11}), $\cos\theta(s)\not=0$. This implies
that $\theta(s)$ is a bounded function with
$-\pi/2<\theta(s)<\pi/2$. If $\bar{s}=\infty$ and as $z(s)>0$, then
both functions $\theta'(s)$ and $z'(s)$ go to $0$ as
$s\rightarrow\infty$. By (\ref{w11}) and (\ref{eq2}),  we have
$(m-1)\cos\theta_1+n=0$ and $\sin\theta_1=0$: contradiction.  This
proves that $\bar{s}<\infty$.

As consequence, $z(s)\rightarrow 0$ since on the contrary,
$\theta'(s)$ would be bounded and $\bar{s}=\infty$. We now use (\ref{integral2}). Letting $s\rightarrow\bar{s}$ and by L'H\^{o}pital rule again, we obtain
$(m-1)\cos\theta_1+n=0$, that is,  $\cos\theta_1\geq -n/(m-1)$.
Finally, $z''(s)=\theta'(s)\cos\theta(s)<0$, that is, $\alpha$ is concave.

\begin{figure}[htbp]\begin{center}
\includegraphics[width=6cm]{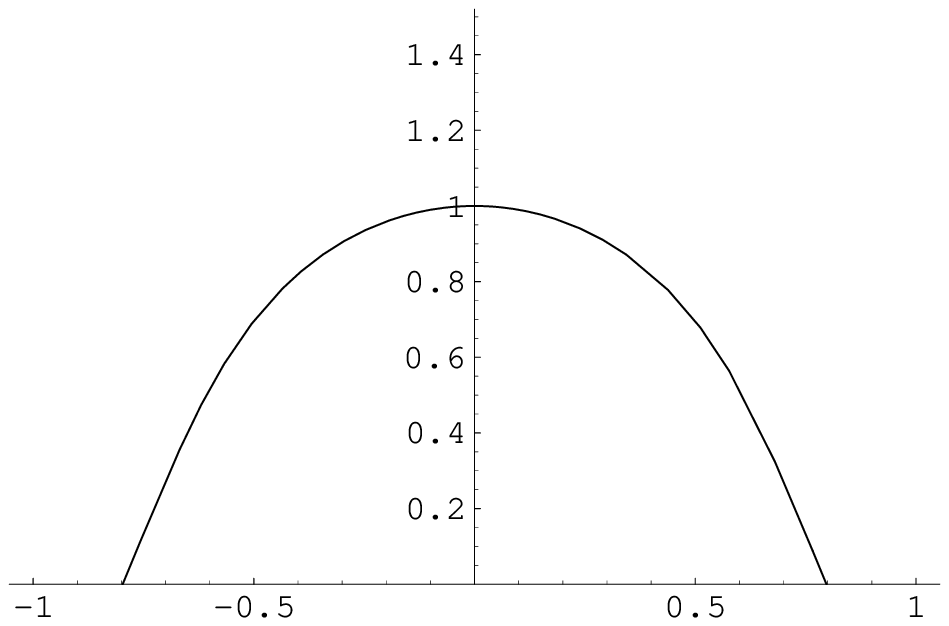}\hspace*{2cm}
\includegraphics[width=4cm]{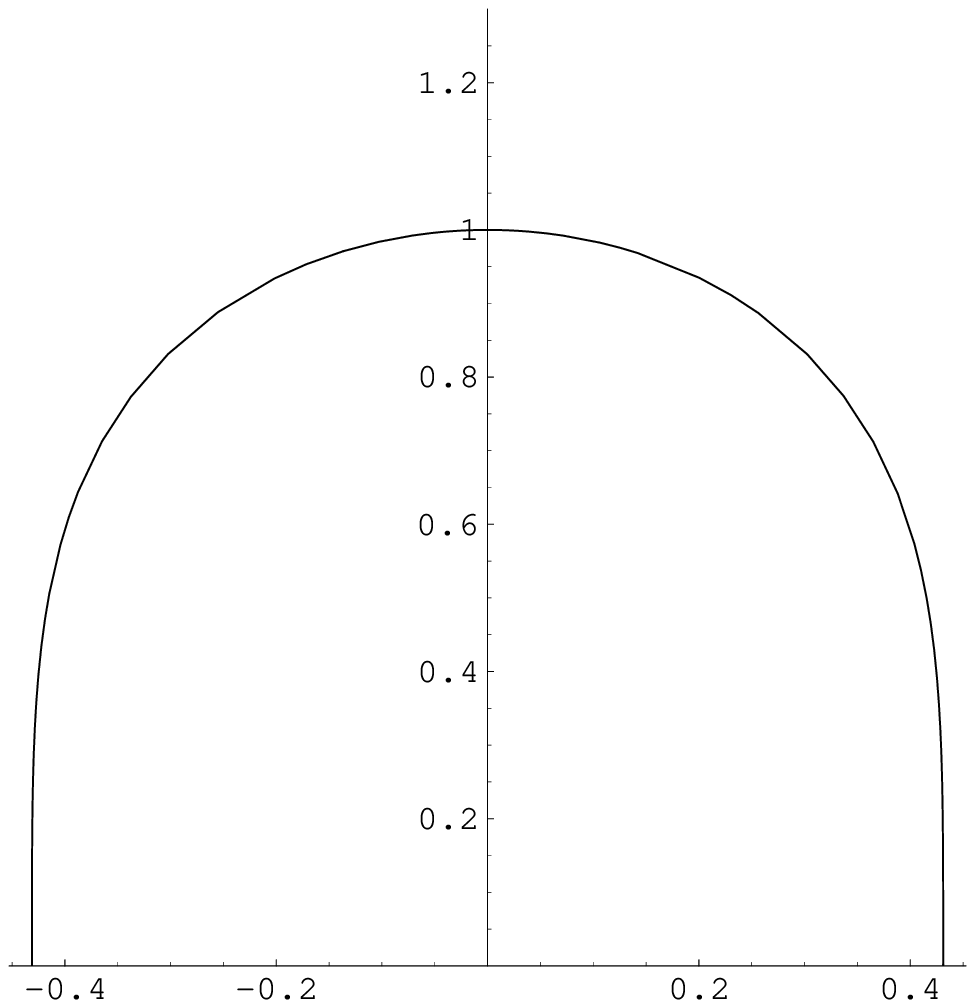}
\end{center}
\begin{center}(a)\hspace*{7cm}(b)\end{center}
\caption{The generating curve of a parabolic surface with
$\kappa_1=m\kappa_2+n$ with $n+m-1<0$. Here $m=-2$, $z_0=1$ and: (a)
$n=1$; (b), $n=0$.}\label{31}
\end{figure}

\begin{theorem} \label{th2} Let $\alpha(s)=(x(s),0,z(s))$ be the generating curve of a
LW-parabolic surface  $S$ in hyperbolic space $\h^3$ whose principal
curvatures  satisfy the relation $\kappa_1=m\kappa_2+n$. Consider
$n\geq 0$ and that  $\theta(0)=0$ in the initial condition
(\ref{eq2}). Assume  $n+m-1<0$.  Then  $\alpha$  is a concave graph
on some bounded interval of $L$ with one maximum and it intersects
$L$ with a contact angle $\theta_1$, $\cos\theta_1= -n/(m-1)$ (see
Fig. \ref{31}, case (a)). In the particular case that $n=0$, then
$\alpha$ meets orthogonally $L$ (see Fig. \ref{31}, case (b)).
 \end{theorem}

As conclusion of this section, we point out the following:

\begin{corollary}\label{mncor} Let $S$ be a parabolic surface in $\h^3$ that
satisfies $\kappa_1= m\kappa_2+n$, for some constants $m$ and $n$.
Consider $\alpha(s)=(x(s),0,z(s))$ the  generating curve. If
$\theta(0)=0$ in the initial condition (\ref{eq2}), then:
\begin{enumerate}
\item The asymptotic boundary $\partial_{\infty} S$ of $S$ is  $\{\infty\}$,
two tangent circles or one circle.
In the latter case, the surface must be umbilical.
\item The surface $S$ is complete.
\item If $S$ is embedded, then it is a graph on $\s^2_{\infty}$.
\end{enumerate}
\end{corollary}


\section{Parabolic surfaces satisfying $a H+bK=c$.}\label{sw2}


In this section we consider parabolic $LW$-surfaces that satisfy the relation
(\ref{wein2}). As in Section \ref{sw1}, we shall consider generating
curves $\alpha$  with some horizontal tangent line. Recall that this
means that   $\theta_0=0$ on the initial condition  (\ref{eq2}). We
also discard the trivial $LW$-surfaces that satisfy equation
(\ref{wein2}), that is, the cases that $a$ or $b$ are zero.
Depending   if the constant $c$, we have two
possibilities:  (i)  $c=0$. Without loss of generality, we assume that
$a=2$;  (ii) $c\not=0$. Then  we take $c=1$.

\subsection{Case $2H+bK=0$}

 Equation (\ref{wein22}) writes as
\begin{equation}\label{w22}
\Big(1+b\cos\theta(s)\Big)z(s)\theta'(s)+2\cos\theta(s)-b\sin^2\theta(s)=0.
\end{equation}
 Let us
denote $(-\bar{s},\bar{s})$ the maximal domain of the solutions of
(\ref{w22}) with initial conditions (\ref{eq2}).
At $s=0$, we have $(1+b)z_0\theta'(0)+2=0$. We know from Lemma \ref{cero}
that $\theta'(s)\not=0$ unless that $\alpha$ is a straight line.
On the other hand, for any solution of (\ref{w22}),
the function $1+b\cos\theta(s)$ can not vanish, since on the contrary, we would have
$\cos\theta(s)=-1/b$ and (\ref{w22}) gives $-b^2-1/b=0$: contradiction.

As $\theta'(0)=-2/(z_0(1+b))$, we have the  subcases $b+1>0$ and $b+1<0$
depending on the sign of $\theta'(0)$. In both settings, let $\theta_1=\lim_{s\rightarrow\bar{s}}\theta(s)$.

\begin{enumerate}
\item Subcase $b>-1$. Then  $\theta(s)$ is a decreasing function.
If $\cos\theta(s)=0$, then (\ref{w22}) yields
$z(s)\theta'(s)-b=0$. Thus, if $b\geq 0$ then $\cos\theta(s)\not=0$ and
$-\pi/2<\theta(s)<\pi/2$.
\begin{enumerate}
\item Let $b\geq 0$. For $s>0$,
$z'(s)=\sin\theta(s)<0$, that is, $z(s)$ is decreasing. If $\bar{s}=\infty$,
 then $\lim_{s\rightarrow\infty}z'(s)=0$.   But (\ref{alpha}) yields $\theta_1=0$: contradiction. Therefore $\bar{s}<\infty$.  If $z(s)\rightarrow z(\bar{s})>0$,
then $\theta'(s)\rightarrow -\infty$ as $s\rightarrow\bar{s}$.
Using (\ref{w22}), $\lim_{s\rightarrow\bar{s}}1+b\cos\theta(s)= 0$,
which it is a contradiction because $1+b\cos\theta(s)\geq 1$.  As conclusion,
 $z(\bar{s})=0$.  From
(\ref{w22}), we obtain that in the contact point between $\alpha$
and the line $L$, both curves make an angle $\theta_1$ such that
$2\cos\theta_1-b\sin^2\theta_1=0$. On the other hand,
$x'(s)=\sin\theta(s)\not=0$ and so $\alpha$ is a graph on  some
bounded interval of $L$. In the particular case that $b=0$, that is,
$S$ is a minimal surface, then  $\alpha$ is a curve that is a graph
on $L$ and it meets $L$ at right angle: this was done in \cite{go}.

\item Consider $-1<b<0$. If $\cos\theta(s)\not=0$, then $\bar{s}<\infty$ as above. We
follow the same reasoning. If $z(\bar{s})>0$ then
$\theta'(s)\rightarrow -\infty$. As
$2\cos\theta(s)-b\sin^2\theta(s)\rightarrow -(1+b^2)/b>0$, then
$\lim_{s\rightarrow\bar{s}}\Big(1+b\cos\theta(s)\Big)=0$. We use the
L'H\^opital rule,
$$\lim_{s\rightarrow\bar{s}}(1+b\cos\theta(s))\theta'(s)=\frac{\lim_{s\rightarrow\bar{s}}
b\theta'(s)\sin\theta(s)}{\lim_{s\rightarrow\bar{s}}\frac{\theta''(s)}{\theta'(s)^2}}>0,$$
since $\theta''(s)<0$ near $\bar{s}$ in contradiction with
(\ref{w22}). Thus $z(\bar{s})=0$ and using (\ref{w22}) again, we
obtain $2\cos\theta_1-b\sin^2\theta_1=0$: contradiction, since
$b<0$.

As conclusion, the function $\theta(s)$ reaches
the value $-\pi/2$ at some point. However, $\theta(s)>-\pi$ using (\ref{w22}) again.
In the case that
 $\bar{s}=\infty$, then   $z'(s)\rightarrow 0$, that is, $\theta_1=-\pi$. But equation
(\ref{w22}) and the fact that $\theta'(s)\rightarrow 0$ gives a
contradiction. Therefore, $\bar{s}<\infty$. We prove that
$z(\bar{s})=0$. On  the contrary, that is, $z(\bar{s})>0$, then
$\theta'(\bar{s})=-\infty$  and for $\theta(s)<-\pi/2$, we would
have
$$\theta'(s)\geq \frac{-2\cos\theta(s)+b\sin^2\theta(s)}{z(s)}\geq
\frac{-2\cos\theta(s)+b\sin^2\theta(s)}{z_0},$$ and  $\theta'(s)$
would be bounded. This contradiction proves the claim on
$z(\bar{s})$. The angle $\theta_1$ which $\alpha$ intersects $L$
satisfies $2\cos\theta_1-b\sin^2\theta_1=0$ by  using (\ref{w22})
again. As
 $x'(s)$ vanishes at some point, then $\alpha$ is not a graph on $L$.
\end{enumerate}

\item Subcase $b<-1$. Now $\theta(s)$ is a strictly increasing function in its
domain. From (\ref{w22}) and as $b<0$,  the function $\cos\theta(s)$
can not vanish. Thus,
  $\theta(s)$ is bounded by $-\pi/2<\theta(s)<\pi/2$.   For $s>0$,
 $z(s)$ is an increasing function. Assuming that $\bar{s}=\infty$,
we will arrive to a contradiction. In such case,
$$\lim_{s\rightarrow\infty}\theta'(s)=0$$
and $\theta''(s)$ is negative  near $s=\infty$.
A differentiation of (\ref{w22}) leads to
\begin{equation}\label{segunda}
(1+b\cos\theta(s))\Big[z(s)\theta''(s)-\sin\theta(s) \theta'(s)\Big]
-b\sin\theta(s) z(s)\theta'(s)^2=0.
\end{equation}
It follows from (\ref{segunda}) that
$\theta''(s)$ is positive  near to $s=0$. Then  $\theta''(s)$ must vanish at some number $s$. However, if $\theta''(s)=0$ it follows from (\ref{segunda}) and the
fact that $1+b\cos\theta(s)<0$ that
for this number $s$, we have
$$0=-\sin\theta(s)\theta'(s)\Big[(1+b\cos\theta(s))+b z(s)\theta'(s)\Big]>0.$$
This contradiction proves  that $\bar{s}<\infty$. This means that
the surface $S$ is not complete. Moreover
$\lim_{s\rightarrow\bar{s}}\theta'(s)=\infty$ and from
(\ref{alpha}), the function $z(s)$ is bounded. Letting
$s\rightarrow\bar{s}$ in  (\ref{w22}) we obtain that
$1+b\cos\theta_1=0$. On the other hand, $x'(s)\not=0$ and so,
$\alpha$ is a graph over a bounded interval of $L$, and as
$z''(s)=\theta'(s)\cos\theta(s)>0$, then $\alpha$ is a convex graph.
\end{enumerate}

\begin{figure}[htbp]\begin{center}
\includegraphics[width=6cm]{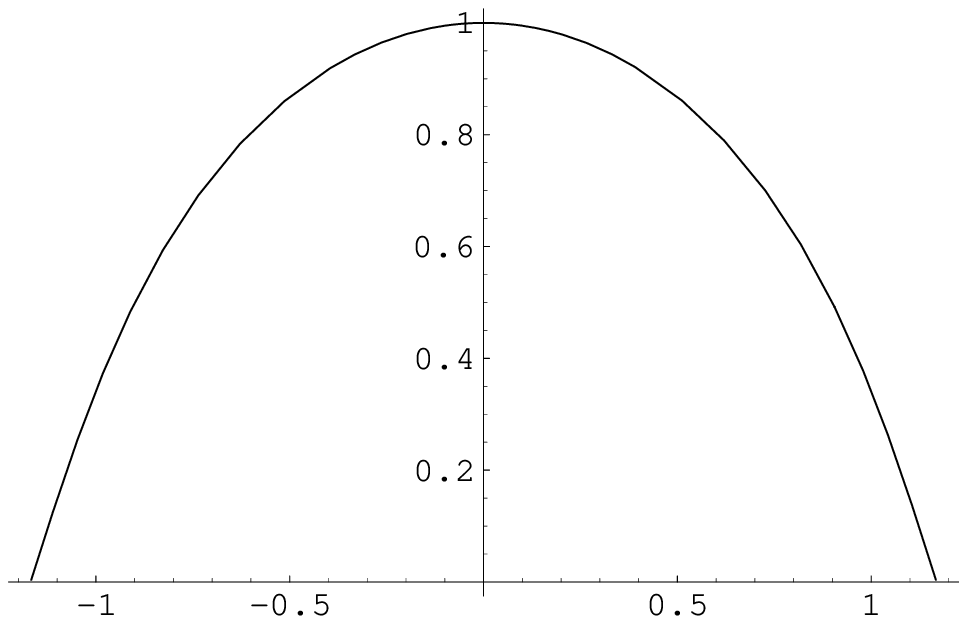}
\hspace*{2cm}\includegraphics[width=6cm]{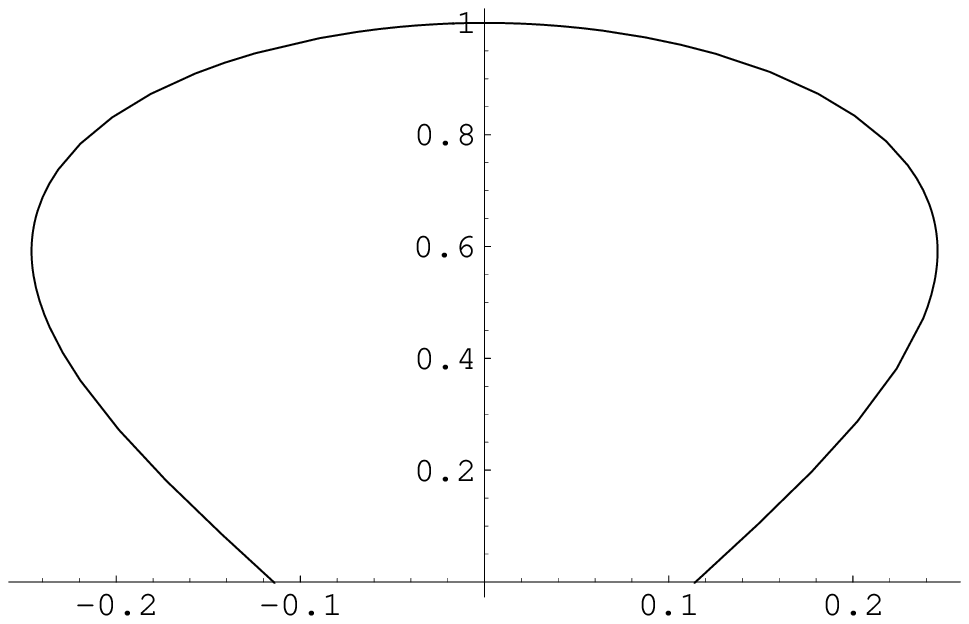}\end{center}
\begin{center}(a)\hspace*{8cm}(b)\end{center}
\caption{The generating curve of a parabolic surface with $2H+bK=0$.
 Here $z_0=1$ and $\theta(0)=0$. In the case
(a), $b=1$ and in the case (b), $b=-0,7$.}\label{51}
\end{figure}

\begin{figure}[htbp]\begin{center}
\includegraphics[width=7cm]{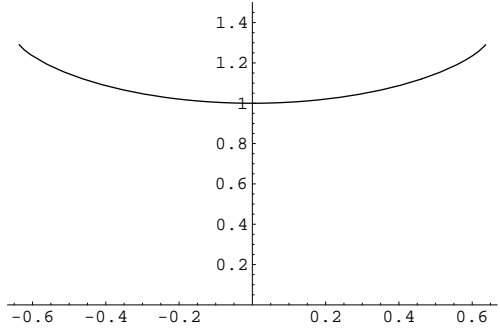}\end{center}
\caption{The generating curve of a parabolic surface with $2H+bK=0$
with $b=-3$.
 Here $z_0=1$ and $\theta(0)=0$. }\label{61}
\end{figure}

\begin{theorem} Let $\alpha(s)=(x(s),0,z(s))$ be the generating curve of a
parabolic surface  $S$ in hyperbolic space $\h^3$ that
satisfies (\ref{alpha})-(\ref{wein22}). Assume that  $\theta(0)=0$. If $2 H+bK=0$ then it
holds the following:
\begin{enumerate}
\item If $b>-1$, $\alpha$ has one maximum and  intersects $L$ at an angle
$\theta_1$ such that $2\cos\theta_1-b\sin^2\theta_1=0$. Moreover, if
$b\geq 0$, $\alpha$ is a concave graph on a bounded interval of $L$,
whereas if  $-1<b<0$, $\alpha$ is not a graph.  See Fig. \ref{51}
cases (a) and (b) respectively.

\item If $b<-1$, then $\alpha$ is a convex graph on a bounded interval of
$L$ with a minimum. At the end points of this interval, $\alpha$
makes an angle $\theta_1$ with the horizontal direction such that
$1+b\cos\theta_1=0$. The surface is not complete. See Fig. \ref{61}.
\end{enumerate}
\end{theorem}

\begin{remark} According to the pictures obtained by the
computer, it seems that in the range $-1<b<0$ there  exists a value
$b_0$
 such  that the curve $\alpha(s)=\alpha(s;b)$ is embedded if $b>b_0$ and it
 is not embedded for $b<b_0$.
However, we have not been able to show the  existence of such number
$b_0$.  In particular, for this value $b_0$,
$\partial_{\infty}\alpha$ is exactly one point, and the
corresponding parabolic surface has one circle as asymptotic
boundary.
\end{remark}

\subsection{Case $aH+bK=1$}

In this subsection we consider parabolic $LW$-surfaces that satisfy
equation (\ref{wein22}) with $c=1$. We discard the trivial
$LW$-surfaces, that is, that $a$ or $b$ are $0$. We also exclude
that situation that for some $s$, $\theta'(s)=0$, and then $\alpha$
would be a straight line by Lemma \ref{cero}. Moreover, and after a
change of orientation on $S$, we can assume that $a>0$.
 Then  equation (\ref{wein22}) writes
\begin{equation}\label{w23}
\left(\frac{a}{2}+b\cos\theta(s)\right)z(s)\theta'(s)+a\cos\theta(s)-b\sin^2\theta(s)=1
\end{equation}
or
\begin{equation}\label{w233}
\theta'(s)=2\frac{1-a\cos\theta(s)+b\sin^2\theta(s)}{z(s)(a+2b\cos\theta(s))}.
\end{equation}

When $a^2+4b^2+4b=0$, we can obtain explicit solutions of (\ref{w23}). Exactly,
equation (\ref{w23}) reduces into
$$-2bz(s)\theta'(s)=a+2b\cos\theta(s).$$
By differentiation with respect to $s$, we obtain
$z(s)\theta''(s)=0$, that is, $\theta'(s)=0$. Since $\theta'(s)$
describes the Euclidean curvature of $\alpha$, we conclude that
$\alpha$ parametrizes a Euclidean circle in the $xz$-plane $P$. This
circle may not to be completely included in the halfspace
${\r}^3_+$.

\begin{theorem} Let $\alpha(s)=(x(s),0,z(s))$ be the generating curve of a
parabolic surface  $S$ in hyperbolic space $\h^3$ that satisfies
(\ref{alpha})-(\ref{wein22}) with $\theta(0)=0$. Assume  $a H+bK=1$
and that $a^2+4b^2+4b=0$. Then $\alpha$ describes an open of a
Euclidean circle in the $xz$-plane $P$. \end{theorem}

We point out that if $\partial_{\infty}\alpha=\emptyset$, that is,
$\alpha$ does not intersect $L$, then the resulting surfaces is
one of the type obtained in Theorem \ref{principal}.  From now, and
along this section, we assume $a^2+4b^2+4b\not=0$.

 Let us
denote $(-\bar{s},\bar{s})$ the maximal domain of the solutions.
At $s=0$, we have
 $$\theta'(0)=\frac{2}{z_0}\frac{1-a}{a+2b}.$$
We study the different settings that appear depending on the sign of
$\theta'(0)$. Moreover, and as $\theta'(s)\not=0$,
 the numerator
can not vanish, that is,
\begin{equation}\label{coseno}
a+2b\cos\theta(s)\not=\pm\sqrt{a^2+4b^2+4b}
\end{equation}
 and it has the same sign
that at $s=0$, that is, as $1-a$. By the monotonicity of
$\theta(s)$, let $\theta_1=\lim_{s\rightarrow\bar{s}}\theta(s)$.
First, we do the computations and next, we state the results. We
will need the computation of the second derivative of $\theta''(s)$:
\begin{equation}\label{dsegunda}
-\theta'(s)\sin\theta(s)\Big[b
\theta'(s)+\left(\frac{a}{2}+b\cos\theta(s)\right)\Big]+
\left(\frac{a}{2}+b\cos\theta(s)\right)z(s)\theta''(s)=0.
\end{equation}
We begin with the case $0<a<1$.

\begin{enumerate}
\item Case $1-a>0$ and $a+2b<0$. Then $\theta'(0)<0$ and $\theta(s)$ is strictly decreasing.
If $\cos\theta(s)=0$ at some point $s$,
then (\ref{w23}) gives $a z(s) \theta'(s)-b-1=0$. Thus, if
$b\geq -1$, $\cos\theta(s)\not=0$ and $-\pi/2<\theta(s)<\pi/2$.
In the case that $b<-1$ and as $a+2b\cos\theta(s)<0$,
 it follows from (\ref{w23}) that  $a\cos\theta(s)-b\sin^2\theta(s)-1<0$ for any value of $s$. In particular, we have        $\cos\theta(s)\not=0$ for any $s$ again. This proves that
$x'(s)=\cos\theta(s)\not=0$ and so,  $\alpha$ is a graph on $L$. This
graph is concave since $z''(s)=\theta'(s)\cos\theta(s)<0$.
Moreover, this implies that $\bar{s}<\infty$ since on the contrary, and as $z(s)$ is
decreasing with $z(s)>0$, we would have $z'(s)\rightarrow 0$,
that is, $\theta(s)\rightarrow 0$: contradiction.

For $s>0$, $z'(s)=\sin\theta(s)<0$ and $z(s)$
 is strictly decreasing. Set $z(s)\rightarrow z(\bar{s})\geq 0$. The two roots of
$4b^2+4b+a^2=0$ on $b$ are $b=-\frac12(1\pm\sqrt{1-a^2})$. Moreover,
and from $a+2b<0$, we have
$$-\frac12(1+\sqrt{1-a^2})<\frac{-a}{2}<-\frac12(1-\sqrt{1-a^2}).$$
\begin{enumerate}
\item Subcase $b< -(1+\sqrt{1-a^2})/2$. With this assumption,
$a^2+4b^2+4b>0$. From (\ref{coseno}) and the fact that $a<1$, we
obtain
\begin{equation} \label{casouno}
a+2b\cos\theta(s)<-\sqrt{a^2+4b^2+4b}.
\end{equation}
If $z(\bar{s})>0$, then
$\lim_{s\rightarrow\bar{s}}\theta'(s)=-\infty$. In particular and
from (\ref{w233}), $a+2b \cos\theta(\bar{s})=0$: contradiction with
(\ref{casouno}). Thus, $z(\bar{s})=0$ and $\alpha$ intersects $L$
with an angle
 $\theta_1$
satisfying $a\cos\theta_1-b\sin^2\theta_1-1=0$.
\item Subcase $-(1+\sqrt{1-a^2})/2< b<-a/2$. Now $a^2+4b^2+4b<0$. The function
 $1-a\cos\theta(s)+b\sin^2\theta(s)$ is strictly decreasing and its value
 at $\bar{s}$ satisfies $\cos\theta(s)>-a/2b$. Thus
 \begin{equation}\label{casodos}
 1-a\cos\theta(s)+b\sin^2\theta(s)\geq \frac{a^2+4b^2+4b}{4b}>0.
 \end{equation}
  Assume $z(\bar{s})=0$.
Then (\ref{casodos}) and (\ref{w233}) imply that
$\theta'(\bar{s})=-\infty$. On the other hand,  using
(\ref{dsegunda}) and (\ref{w233}), we have
 $$\frac{\theta''(s)}{\theta'(s)^2}
 =\frac{b\sin\theta(s)}{z(s)\left(\frac{a}{2}+b\cos\theta(s)\right)}+
\frac{\sin\theta(s)\left(\frac{a}{2}+b\cos\theta(s)\right)}{1-a\cos\theta(s)+b\sin^2\theta(s)}.
$$
From this equation and as  $\sin\theta(\bar{s})\not=0$, we conclude
$$\lim_{s\rightarrow\bar{s}}\frac{\theta''(s)}{\theta'(s)^2}=-\infty.$$
On the other hand, using L'H\^opital rule, we have
$$\lim_{s\rightarrow\bar{s}}z(s)\theta'(s)=
\lim_{s\rightarrow\bar{s}}-\frac{\sin\theta(s)}{\frac{\theta''(s)}{\theta'(s)^2}}=0.$$
As the numerator in  (\ref{w233}) is bounded from below for a
positive number, see (\ref{casodos}), we obtain a contradiction by
letting $s\rightarrow\bar{s}$.
 Thus, $z(\bar{s})>0$. This means
that $\lim_{s\rightarrow\bar{s}}\theta'(s)=-\infty$ and from
(\ref{w233}), that
$$\lim_{s\rightarrow\bar{s}}\left(\frac{a}{2}+b\cos\theta(s)\right)=0.$$

\end{enumerate}
\item Case $1-a>0$ and $a+2b>0$. Then $\theta'(0)>0$ and $\theta(s)$ is strictly increasing.
We distinguish  two possibilities:
\begin{enumerate}
\item Subcase $a-2b>0$. We prove that $\theta(s)$ reaches the value $\pi$. On the
contrary, $\theta(s)<\pi$ and $z(s)$ is an increasing function. The
hypothesis $a-2b>0$ together $a+2b>0$ implies that
$a+2b\cos\theta(s)\geq \delta>0$ for some number $\delta$. From
(\ref{w233}), $\theta'(s)$ is bounded and then $\bar{s}=\infty$. In
particular, $\lim_{s\rightarrow\infty}\theta'(s)=0$. As both $a-2b$
and $a+2b$ are negative numbers, the function $b\theta'(s)+(a+2b\cos
\theta(s))$ is positive near $\bar{s}=\infty$. Then using
(\ref{dsegunda}), from a certain big value of $s$, $\theta''(s)$ is
positive, which it is impossible. As conclusion, $\theta(s)$ reaches
the value $\pi$ at some $s=s_0$. By Lemma \ref{si1},  $\alpha$ is
symmetric with respect to the line $x=x(s_0)$ and the velocity
vector of $\alpha$ rotates until to the initial position. This means
that $\alpha$ is invariant by a group of horizontal translations.
\item Subcase $a-2b\leq 0$. As $\theta'(s)>0$, Equation (\ref{w233}) says that
$\cos\theta(s)\not=-1$, and so, $\theta(s)$ is bounded by
$-\pi<\theta(s)<\pi$.  As in the above subcase, if $\bar{s}=\infty$,
then $\theta'(s)\rightarrow 0$, and this is a contradiction. Then
$\bar{s}<\infty$ and $\lim_{s\rightarrow\bar{s}}\theta'(s)=\infty$.
In particular, $\cos\theta(\bar{s})=-a/(2b)$ and $\theta(s)$ reaches
the value $\pi/2$.
\end{enumerate}
\end{enumerate}

\begin{figure}[htbp]\begin{center}
\includegraphics[width=6cm]{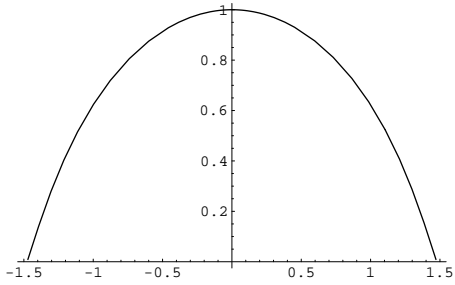}
\hspace*{2cm}\includegraphics[width=6cm]{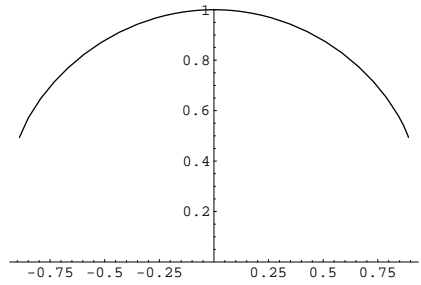}\end{center}
\begin{center}(a)\hspace*{8cm}(b)\end{center}
\caption{The generating curve of a parabolic surface with $aH+bK=1$,
with $0<a<1$ and $a+2b<0$.
 Here $z_0=1$, $\theta(0)=0$ and $a=0.5$. In the case
(a), $b=-1$ and in the case (b), $b=-0.8$.}\label{71}
\end{figure}

\begin{figure}[htbp]\begin{center}
\includegraphics[width=6cm]{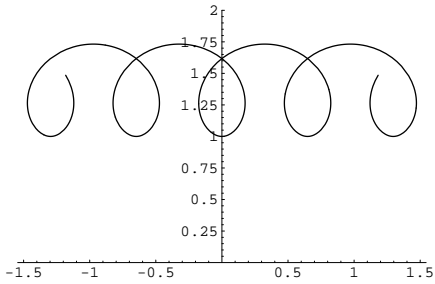}
\hspace*{2cm}\includegraphics[width=6cm]{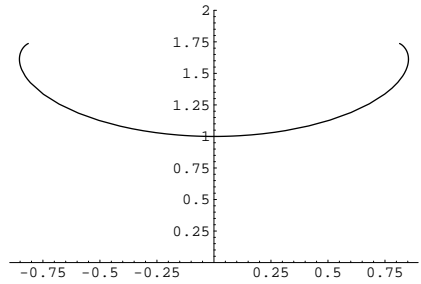}\end{center}
\begin{center}(a)\hspace*{8cm}(b)\end{center}
\caption{The generating curve of a parabolic surface with $aH+bK=1$,
with $0<a<1$ and $a+2b>0$.
 Here $z_0=1$, $\theta(0)=0$ and $a=0.5$. In the case
(a), $b=-0.2$ and in the case (b), $b=0.3$.}\label{81}
\end{figure}

\begin{theorem} Let $\alpha(s)=(x(s),0,z(s))$ be the generating curve of a
parabolic surface  $S$ in hyperbolic space $\h^3$ that satisfies
(\ref{alpha})-(\ref{wein22}) with $\theta(0)=0$. Assume  $a H+bK=1$
and that $0<a<1$.
\begin{enumerate}
\item If $a+2b<0$, $\alpha$ has one maximum and $\alpha$ is a
concave graph in some bounded interval of $L$.  If $b<
-(1+\sqrt{1-a^2})/2$, $\alpha$ intersects $L$
 at an angle $\theta_1$ such that
$2\cos\theta_1-b\sin^2\theta_1=0$. If $-(1+\sqrt{1-a^2})/2<b<-a/2$,
then $\alpha$ does not intersect $L$.  In the latter case, the
surface is not complete. See Fig. \ref{71} cases (a) and (b)
respectively.

\item Assume $a+2b>0$. If $a-2b>0$, then $\alpha$ is invariant by a
group of translations in the $x$-direction, $\alpha$ has
self-intersections and it presents one maximum and one minimum in
each period. The velocity $\alpha'$ turns around the origin. If
$a-2b\leq 0$, then $\alpha$ has a minimum and it is not a graph on
$L$. In this case, the surface $S$ is not complete. See Fig.
\ref{81}, cases (a) and (b) respectively.
\end{enumerate}
\end{theorem}

If   $a=1$ and since $\theta'(0)=0$, we obtain that  $\alpha$ is a
straight line. Thus $\alpha$ is horizontal and $S$ is a horosphere.

\begin{theorem} Let $\alpha(s)=(x(s),0,z(s))$ be the generating curve of a
parabolic surface  $S$ in hyperbolic space $\h^3$ that satisfies
(\ref{alpha})-(\ref{wein22}) with $\theta(0)=0$. Assume  $ H+bK=1$.
Then $\alpha$ is a horizontal straight line and $S$ is a horosphere.
\end{theorem}

 Consider the setting that $a>1$.

\begin{enumerate}

\item Case $a>1$ and $a+2b<0$. Now,  the function
$a+2b\cos\theta(s)<0$, in particular, $\cos\theta(s)\not=0$. This
means  $-\pi/2<\theta(s)<\pi/2$ and $\alpha$ is a graph on $L$.
Moreover, $\cos\theta(s)>-a/2b$, that is,  $\theta_1<\pi/2$. As
$\theta'(0)>0$, then $\theta(s)$ is a strictly increasing function
and the same occurs for $z(s)$ for $s>0$.  We claim that
$\bar{s}<\infty$. Assuming the contrary and as $\theta(s)<\theta_1$,
we have that $\lim_{s\rightarrow \infty}\theta'(s)=0$. From
(\ref{dsegunda}), $\theta''(s)>0$ in a neighborhood of $\infty$,
which it is impossible. As conclusion, $\bar{s}<\infty$. Then
$\lim_{s\rightarrow\bar{s}}\theta'(s)=\infty$. Since $z(s)$ is
defined in a bounded interval and its derivative is bounded, then
$a/2+b\cos\theta(s)\rightarrow 0$ as $s\rightarrow\bar{s}$.

\item Case $a>1$ and $a+2b>0$. Now $\theta'(s)<0$,
$\theta(s)$ is a decreasing function and  $a+2b \cos\theta(s)>0$.
From (\ref{w23}), if $a-2b>0$, then $\cos\theta(s)\not=-1$ for any
$s$ and if $a-2b\leq 0$, then $\cos\theta(s)\not=0$. As conclusion,
$\theta(s)$ is a bounded function. We prove that $\bar{s}<\infty$.
On the contrary and since $z(s)>0$,  both functions $z'(s)$ and
$\theta'(s)$ go to $0$ as $s\rightarrow \infty$. This means that
$\theta_1=-\pi$ and from (\ref{w23}) and letting
$s\rightarrow\infty$ we conclude
 that $a=1$. This contradiction proves that
$\bar{s}<\infty$. We have two possibilities: either
$\theta'(\bar{s})=-\infty $ or $z(\bar{s})=0$. The first case is
impossible because $\phi(s):=a+2b\cos\theta(s)\geq\delta>0$, for
some number $\delta>0$: if $a-2b>0$, then $\phi(s)\geq a-2b$ and if
$a-2b\leq 0$, then $\cos\theta(s)>0$ and so, $\phi(s)\geq a$. As
conclusion, $z(\bar{s})=0$, that is, $\alpha$ intersects the line
$L$.
\end{enumerate}

\begin{figure}[htbp]\begin{center}
\includegraphics[width=6cm]{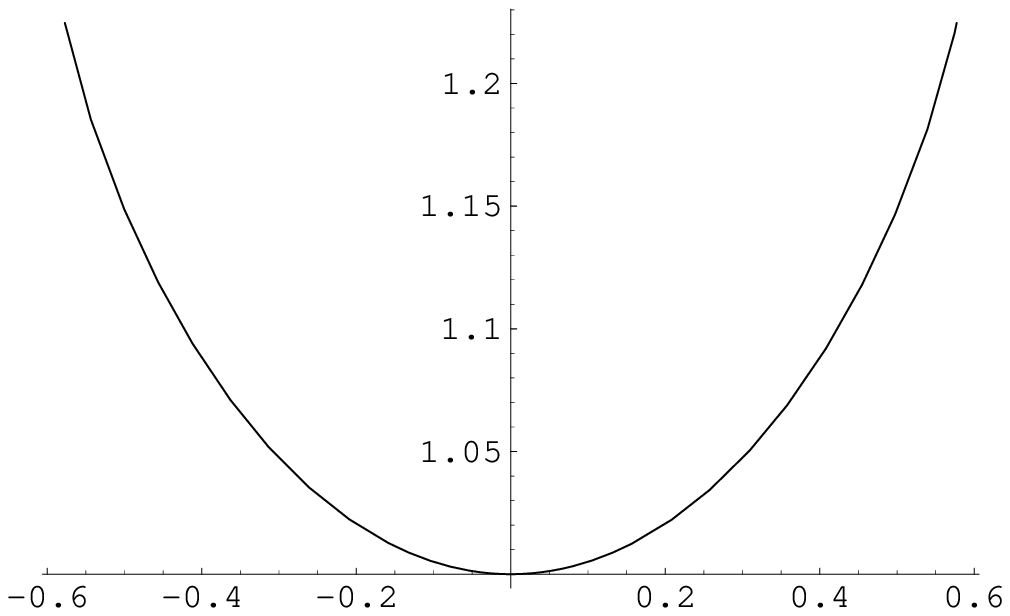}
\hspace*{2cm}\includegraphics[width=6cm]{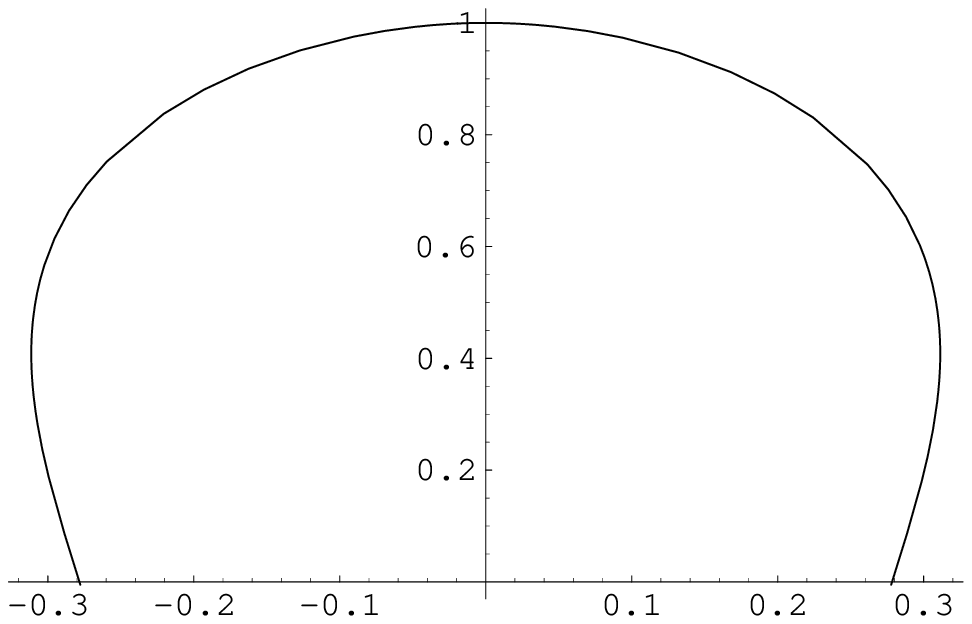}\end{center}
\begin{center}(a)\hspace*{8cm}(b)\end{center}
\caption{The generating curve of a parabolic surface with $aH+bK=1$,
with $a>1$.
 Here $z_0=1$ and $\theta(0)=0$. In the case
(a), $a=2$ and $b=-2$ and in the case (b), $a=4$ and
$b=-1.5$.}\label{91}
\end{figure}

\begin{theorem} Let $\alpha(s)=(x(s),0,z(s))$ be the generating curve of a
parabolic surface  $S$ in hyperbolic space $\h^3$ that satisfies
(\ref{alpha})-(\ref{wein22}) with $\theta(0)=0$. Assume  $a H+bK=1$
and that $a>1$.
\begin{enumerate}
\item If $a+2b<0$, $\alpha$ has one minimum and it is a convex graph in some
bounded interval of $L$. See Fig. \ref{91} case (a). The surface is
not complete.

\item If $a+2b>0$, then $\alpha$ has a maximum and intersects $L$ with an angle
$\theta_1$ such that $1-a\cos\theta_1+b\sin^2\theta_1=0$. See Fig.
\ref{91}, case (b).
\end{enumerate}
\end{theorem}

\section{Final remarks}\label{final}

As we pointed out in Introduction of this work, we have considered
parabolic surfaces in $\h^3$ whose generating curve has a tangent
line parallel to the plane $\Pi$. This means that
$\theta(0)=0$ in the initial condition (\ref{eq2}). However, similar
reasonings can done in the general case of $\theta_0$. In this
section, we describe briefly what occurs for parabolic surfaces that
satisfy the Weingarten condition $\kappa_1=m\kappa_2+n$. We left the
reader the setting of $aH+bK=c$. We return in Section \ref{sw1}
with each one of the cases discussed there. We  have to see if the
angle function $\theta(s)$ takes all the possible values in the
range $[0,2\pi]$. The reasonings are similar, and we only sketch
them.

\begin{figure}[htbp]\begin{center}
\includegraphics[width=6cm]{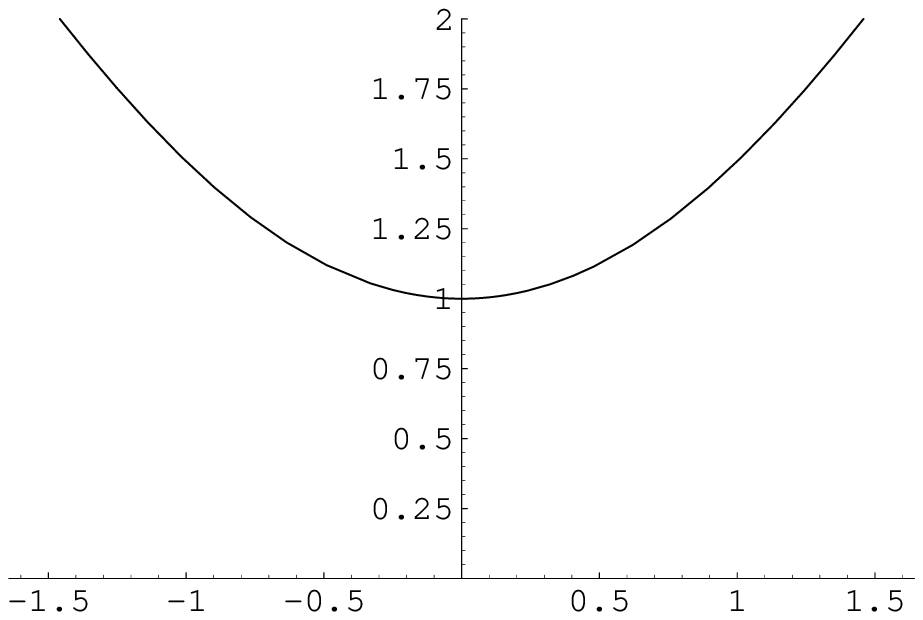}
\hspace*{2cm}\includegraphics[width=6cm]{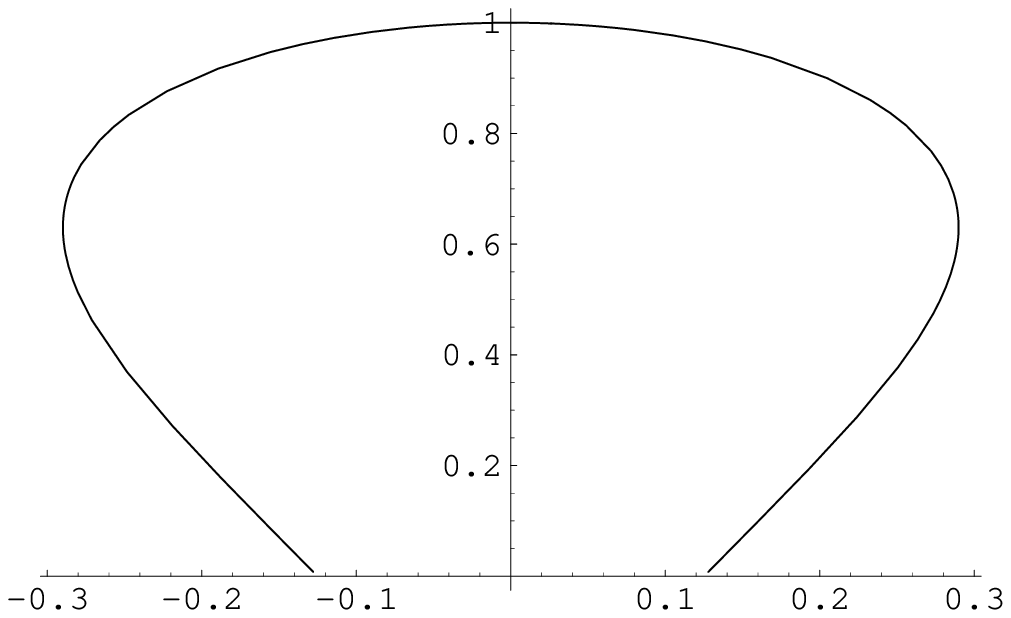}\end{center}
\begin{center}(a)\hspace*{8cm}(b)\end{center}
\caption{The generating curve of a parabolic surface with
$\kappa_1=m\kappa_2+n$.  Here $z_0=1$ and $\theta(0)=\pi$.   In the
case (a), $m=3$, $n=1$; in the case (b), $m=-2$, $n=1$.}\label{10}
\end{figure}

\begin{enumerate}
\item Case $n+m-1>0$. In the subcase $m<n+1$, the function $\theta(s)$ takes all
the possible values of the interval $[0,2\pi]$. Then, let us
consider $m\geq n+1>0$. First, we assume $n=0$. Let
$\theta_0=\pi/2$. Then the solution of (\ref{alpha})-(\ref{eq2}) is
$\alpha(s)=(0,0,s+z_0)$, that is, $\alpha$ is a vertical line and
$S$ is a geodesic plane. On the other hand, if $\theta_0=\pi$ and
$(x(s),z(s),\theta(s))$ is the corresponding solution, then
$(x(-s),z(-s),\theta(-s)-\pi)$ is the solution for $\theta_0=0$, and
this case has been  already  studied.

Let $n>0$. From Theorem \ref{th1}, we only have to consider
$\theta_0=\pi$. This case is similar to the considered one when
$\theta_0=0$. As $\theta'(0)<0$, $\theta(s)$ is a decreasing
function and $\cos\theta(s)<-n/(m-1)$. In particular, $x'(s)\not=0$
and $z''(s)>0$. This implies that $\alpha$ is a convex graph on $L$
and  at infinity, $\theta(s)$ goes to the value
$\theta_1\in(\pi/2,\pi)$ with $\cos\theta_1=n/(m-1)$. See Fig.
\ref{10}, (a).

\item Case $n+m-1=0$. As $n\geq 0$, then $m-1<0$
(recall that we discard the umbilical case $m=1, n=0$). Assume
$\theta_0=\pi$. Then by (\ref{w11}), $\theta'(s)>0$ and
$\cos\theta(s)\not=1$. This means that $0<\theta(s)<2\pi$. For
$s>0$, $z'(s)<0$ and $z(s)\leq z_0$. Thus $\theta'(s)$ is a bounded
function and $\bar{s}=\infty$. Then $\theta'(s)\rightarrow 0$ as
$s\rightarrow\infty$. From (\ref{w11}), this yields that
$$\lim_{s\rightarrow\infty}\cos\theta(s)=0,\hspace*{1cm}
\lim_{s\rightarrow\infty}z(s)=0.$$ Since $\theta(s)$ takes all the
values of $(0,2\pi)$, we have proved:

\begin{theorem}\label{ultimo}
Let $\alpha(s)=(x(s),0,z(s))$ be the generating curve of a parabolic
surface  $S$ in hyperbolic space $\h^3$. Assume that the principal
curvatures of $S$ satisfy the relation $\kappa_1=m\kappa_2+n$ with
$n+m-1=0$. If $\theta(0)\in (0,2\pi)$ in the initial condition
(\ref{eq2}), then
\begin{enumerate}
\item $\alpha$ has one maximum.
 \item $\alpha$ has self-intersections.
 \item $\alpha$ is asymptotic to the line $L$ at infinity, that is,
 $$\lim_{s\rightarrow\pm\infty}z(s)=0.$$
 \end{enumerate}
See Fig. \ref{41}, case (b).
\end{theorem}

\item Case $n+m-1<0$. Recall that Theorem \ref{th2} says that when $\theta_0=0$,
$\theta(s)\in(-\theta_1,\theta_1)$ with $\cos\theta_1=-n/(m-1)$.
Consider then $\theta_0=\pi$. Then $\theta'(0)>0$ and so,
$\theta(s)$ is strictly increasing, with $-1\leq
\cos\theta(s)<-n/(m-1)$. As in Theorem \ref{th2}, one can show that
$\theta(s)$ takes all values of the interval
$(\theta_1,2\pi-\theta_1)$ and that intersects the ideal boundary at
$L$ at two points. Thus the assumption $\theta_0=\pi$ covers all
possibilities of the initial angle. As the angle function $
\theta(s)$ reaches the values $\pi/2$ and $3\pi/2$, $\alpha$ is not
a graph. See Fig. \ref{10}, (b). In the case that $n=0$, we point
out that if $\theta_0=\pi/2$, then $\alpha$ is a vertical line.
\end{enumerate}

\begin{remark} We compare the results obtained in this Section with
 the ones  in Corollary \ref{mncor}. Here, we
have showed the existence of  parabolic surfaces in $\h^3$ that
satisfy the relation $a\kappa_1+b\kappa_2=c$ with the property that
either i) are not complete or ii) are embedded,  complete but not
graphs.
\end{remark}

\end{document}